\documentstyle[12pt]{article}
\newcommand{\beq}{\begin{equation}\label}
\newcommand{\eeq}{\end{equation}}
\newcommand{\bea}{\begin{eqnarray}\label}
\newcommand{\eea}{\end{eqnarray}}
\newcommand{\ar}{\begin{array}}
\newcommand{\er}{\end{array}}

\newcommand{\up}{\stackrel}
\newcommand{\ti}{\widetilde}
\newcommand{\rar}{\rightarrow}

\newcommand{\ov}{\overline}
\newcommand{\wh}{\widehat}
\newcommand{\red}{\hspace*{6mm}}

\newcommand{\al}{\alpha}

\newcommand{\ga}{\gamma}

\newcommand{\de}{\Delta}
\newcommand{\dl}{\delta}

\newcommand{\qed}{\hfill {$\Box$}}

\newcommand{\Lo}[3]{\stackrel{#1}{L}\hspace*{-4mm}\phantom{L}^{#2}_{#3}}
\newcommand{\Uo}[2]{\stackrel{#1}{U}\hspace*{-4mm}\phantom{U}^{#2}}
\newcommand{\R}[3]{\stackrel{#1}{R}\hspace*{-4mm}\phantom{R}^{#2}_{#3}}
\newcommand{\Z}[2]{\stackrel{#1}{Z}\hspace*{-4mm}\phantom{Z}^{#2}}
\newcommand{\W}[2]{\stackrel{#1}{W}\hspace*{-5mm}\phantom{W}^{#2}}
\newcommand{\Uw}[2]{\stackrel{#1}{\wh{U}}\hspace*{-4mm}\phantom{U}^{#2}}
\newcommand{\Ww}[2]{\stackrel{#1}{\wh{W}}\hspace*{-5mm}\phantom{W}^{#2}}
\newcommand{\uo}[2]{\stackrel{#1}{U_0}\hspace*{-5.5mm}\phantom{U}^{#2}}
\newcommand{\wo}[2]{\stackrel{#1}{W_0}\hspace*{-7mm}\phantom{W}^{#2}}
\newcommand{\lo}[3]{\stackrel{#1}{l}\hspace*{-2mm}\phantom{l}^{#2}_{#3}}
\newcommand{\ho}[2]{\stackrel{#1}{\chi}\hspace*{-3mm}\phantom{\chi}^{#2}}
\newcommand{\Po}[3]{\stackrel{#1}{P}\hspace*{-4mm}\phantom{P}^{#2}_{#3}}
\newcommand{\F}[2]{\stackrel{#1}{F}\hspace*{-4mm}\phantom{P}^{#2}}
\newcommand{\Uf}[3]{\stackrel{#1}{U}\hspace*{-4mm}\phantom{U}^{#2}_{#3}}
\newcommand{\Wf}[3]{\stackrel{#1}{W}\hspace*{-5mm}\phantom{U}^{#2}_{#3}}
\newcommand{\ve}[3]{\stackrel{#2}{|#1\rangle}\hspace*{-4mm}\phantom{#1}^{#3}}

\def\nn{\nonumber}

\newtheorem{prop}{Proposition}

\begin{document}
\hyphenation{ope-ra-tor ope-ra-tors ge-ne-ra-tor ge-ne-ra-tors
ope-ra-tion ope-ra-tions}

\thispagestyle{empty}
\vspace*{0.3cm}

\hfill {\begin{tabular}{c}  q-alg/9512030 \\   \\
December 1995 \end{tabular} }

\vspace*{2.5cm}  \begin{center}
{\Large \bf Tensor operators in $R$-matrix approach.}  \end{center}
\vspace{1cm}

\begin{center}  {\bf A. G. Bytsko} $^*$ $^{\#}$  \end{center}
\vspace*{0.1cm}

\begin{center} \begin{tabular}{c}
II. Institut f\"ur Theoretische Physik, \\ Universit\"at Hamburg, \\
Luruper Chaussee 149, \\ 22761 Hamburg, Germany
\end{tabular} \end{center}
 \vspace*{0.4cm}

 \vspace*{5.5cm}

\begin{center} \bf{ Abstract  } \end{center}
\red
The definitions and some properties (e.g. the Wigner-Eckart theorem, the
fusion procedure) of covariant and contravariant $q$-tensor operators for
quasitriangular quantum Lie algebras are formulated in the $R$-matrix
language. The case of $U_q(sl(n))$ (in particular, for $n=2$) is discussed
in more detail.

\vfill

\begin{flushleft}   \rule{7 cm}{0.05 mm} \\
\begin{tabular}{cl}
$^{*}$ & \footnotesize e-mail : \ bytsko@x4u2.desy.de  \\
$^{\#}$ & \footnotesize On leave of absence from \\
& \footnotesize Steklov Mathematical Institute, Fontanka 27, St.-Petersburg,
191011, Russia; \\
& \footnotesize e-mail : \ bytsko@pdmi.ras.ru  \end{tabular}
\end{flushleft}

\newpage
\red {\large\bf 1. Introduction}
\vspace*{2mm}

Quantum groups and quantum algebras which originally appeared in the studies
of integrable systems present nowadays a well developed part of mathematics.
The $R$-matrix method is the most powerful and very elegant tool of this
theory.  On the other hand, the theory of tensor operators which has its
origins in quantum mechanics is also a well-known mathematical construction.
Its modification -- the theory of $q$-tensor operators has been intensively
developing during several last years.

The purpose of this paper is to give a formulation of the theory of $q$-tensor
operators using the $R$-matrix language.

The paper is organized as follows. First, we remind some basic facts
concerning $q$-tensor operators and quantum algebras. Next, we give a
definition of tensor operators in the $R$-matrix language. And further we
discuss from this point of view different aspects of general theory of tensor
operators (Wigner-Eckart theorem, fusion procedure, construction of scalars
and invariants, classical limit).
\vspace*{2mm}

\red {\bf a) Tensor operators and their $q$-analogues.}
\vspace*{2mm}

The most of physical models possess groups of symmetry. To
classify the objects described by such a model (fields, currents etc.)
with respect to action of the group one needs to introduce the notion of
{\it tensor operator} \cite{0}. Let us recall the definition and some
properties of tensor operators.

Let $G$ be a Lie group and $\rho^I$ be some finite dimensional
representation with carrier vector space $V^I$. Then the set of operators
$\{T_m^I\}_{m=1}^{{\rm dim}\,\rho^I}$ is called a tensor operator
(transforming according to the representation $\rho^I$) if
\beq{1}
  U(g)\,T^I_m\,U^{-1}(g) = \sum\limits_n \,T^I_n\,\left(\rho^I(g)\right)_{nm},
\;\;\;\; {\rm for \ all \ } g\in G,
\eeq
where $U(g)$ is a unitary representation of the group $G$ in Hilbert space
$\cal H$.

The definition (\ref{1}) rewritten in terms of the corresponding Lie
algebra $\cal J$ looks as follows
\beq{3}
{} [\,\xi\, ,\, T^I_m \, ] = \sum\limits_n \,T^I_n \, \left(\rho^I(\xi)
\right)_{nm}, \;\;\;\; {\rm for \ all \ } \xi\in{\cal J},
\eeq
where $\rho^I(\xi)\!\equiv\!\frac{1}{i}\frac{d}{d\theta}\rho^I
(\exp\{i\theta\xi\}){\Bigm|}_{\theta=0}$.

{\it Remark}. In what follows we are considering tensor operators $\{T^I_m\}$
as operators acting on the {\it model space} \cite{3}, which is the direct sum
of all irreducible representations of the group $G$ with multiplicity one:
\beq{M}
{\cal M}\, =\, \bigoplus\limits_{I} \; {\cal H}_{I} .
\eeq
Notice that elements of the algebra ${\cal J}$ itself also are operators
acting on ${\cal M}$.

The main property of tensor operators is described by the Wigner-Eckart
theorem \cite{0}, which allows to reduce the problem of finding the
Clebsch-Gordan coefficients (CGC) to the problem of constructing of the
corresponding tensor operators:

{\bf Theorem} (Wigner-Eckart). \ {\it Let
$\{ |I\,m\rangle \}_{m=1}^{{\rm dim}\,\rho^{I}}$ denote the set
of orthonormal basic vectors in the irreducible subspace ${\cal H}_{I}$ of the
model space $\cal M$ and let $\{T_m^J\}_{m=1}^{{\rm dim}\,\rho^{J}}$ be a
tensor operator. Then matrix elements of the components $T_m^{J}$ are
\beq{WE}
 \langle K m^{\prime\prime} |T^{J\;}_{\;m}| {I} m^{\prime} \rangle
= ||T^J||_{I\,K}\; \Bigl\{\!\ar{ccc} I &\!J\!&\!K\\ m^{\prime}\!&\! m\!&\!
m^{\prime\prime}\er\!\Bigr\},
\eeq
where $\Bigl\{\!\ar{ccc} I &\!J\!&\!K\\ m^{\prime}\!&\! m\!&\!
m^{\prime\prime}\er\!\Bigr\}$ are the Clebsch-Gordan coefficients
corresponding to the decomposition of the tensor product
$\rho^J\!\otimes\!\rho^I$ and the coefficients $||T^J||_{I\,K}$ (called
reduced matrix elements) do not depend on $m$, $m^{\prime}$,
$m^{\prime\prime}$.}  \\

The theory of tensor operators can be extended to the case of quantum
Lie algebras \cite{BMT,Bi,MS}. For this purpose the most natural way
is to use the Hopf algebra attributes (i.e., the co-product $\de$, the
antipode $S$ and the co-unit $\epsilon$) of ${\cal J}_q$.

Let $\xi$ be an element of ${\cal J}_q$ with the co-product
$\de\xi=\sum\limits_k {\xi}_k^1 \otimes{\xi}_k^2$ and the antipode $S(\xi)$.
Let $\cal H$ be a certain Hilbert space such that
${\cal J}_q\subset{\rm End}\,{\cal H}$.
The definition of the adjoint action of the algebra reads \cite{BMT,Bi,MS}
\beq{5}
  ({\rm ad}_q \xi)\;\eta = \sum\limits_k \; {\xi}_k^1 \; \eta \;
S({\xi}_k^2), \;\;\;\;\;  \xi\in{\cal J}_q,\;\; \eta\in {\rm End}\,{\cal H} .
\eeq

The analogue of definition (\ref{3}) for $q$-tensor operators looks like
\beq{6}
 ({\rm ad}_q \xi)\, T_m^I
= \sum\limits_n \,T^I_n \, \left(\rho^I(\xi)\right)_{nm},
\;\;\;\; {\rm for \ all \ } \xi\in {\cal J}_q.
\eeq
Here $\rho$ is a finite dimensional representation of ${\cal J}_q$.

{\it Remark}. Let us underline that in (\ref{6}) the elements of the quantum
algebra and the components of tensor operator are regarded as operators
acting on the same Hilbert space $\cal H$ (say, on the model space $\cal M$).

{\it Remark}. The definition (\ref{6}) admits the following vector form
($t$ denotes a transposition)
\beq{2.4a}
\left({\rm ad}_q\,\xi\right)\,\vec{T}^I\, = \, \Bigl(\rho^I(\xi)\Bigr)^{t}
\,\vec{T}^I ,\;\;\;\; \vec{T}^I\equiv \left(\! \ar{c} \scriptstyle T_1^I \\
\vdots \\ \scriptstyle T_N^I \er \!\right) , \;\;\;\; N\equiv {\rm dim}\,
\rho^I,
\eeq
or, equivalently,
\beq{2.4b}
\left({\rm ad}_q\,\xi\right)\,(\vec{T}^I)^t\, =  (\vec{T}^I)^t \,\rho^I(\xi),
\;\;\;\; (\vec{T}^I)^t\equiv ({\scriptstyle T^I_1},..,{\scriptstyle T^I_N}).
\eeq

{\it Remark}. The adjoint action (\ref{5}) has the property
$\left({\rm ad}_q (\xi\,{\xi}^{\prime})\right)\,\eta= ({\rm ad}_q \xi)\,
\left(({\rm ad}_q \xi^{\prime})\,\eta\right)$. Applying it to the
definition (\ref{6}), one gets $\left({\rm ad}_q (\xi\,{\xi}^{\prime})\right)
\,T_m^I=\sum\limits_n T_n^I\,\left(\rho^I(\xi)\,\rho^I({\xi}^{\prime})
\right)_{nm}$. This implies that it is sufficient to consider the definition
(\ref{6}) only for the basic elements of the algebra (say, for the generators
corresponding to simple roots).

{\it Remark}. In the classical limit ($q\rar 1$) the definition (\ref{6})
gives (\ref{3}) since for the generators of classical algebra
$\de\,\xi\!=\!\xi\!\otimes\!1+1\!\otimes\!\xi$, $S(\xi)\!=\!-\xi$
(and due to the previous Remark, it is sufficient to check the coincidence
only on the set of generators).

Similarly to the classical case, matrix elements of the components of the
$q$-tensor operators are proportional to deformed CGC (more on the
generalization of the Wigner-Eckart theorem can be found e.g. in \cite{Bi}).

In general, physical models involve besides tensor operators other objects
which transform as {\it conjugated} tensor operators (these two types of
objects are also called co- and contravariant operators).
Recall that for the classical Lie group $G$ the conjugated tensor operator
$\{{\ov{T}}_m^{\,I}\}_{m=1}^{{\rm dim}\,\rho^I}$ is defined as follows
\cite{0}
\beq{37}
 U(g)\,\ov{T}_m^{\,I}\,U^{-1}(g) =
\sum\limits_n \left({\rho^I}(g^{-1})\right)_{mn} \,\ov{T}^{\,I}_n,
\eeq
or, equivalently, in terms of generators of the corresponding Lie algebra
${\cal J}$ as follows
\beq{co1}
{} [\,\xi\, ,\, \ov{T}_m^{\,I} \, ] = -\sum\limits_n \left(\rho^I(\xi)
\right)_{mn}\, \ov{T}_n^{\,I}.
\eeq
The $q$-analogue of the definition (\ref{co1}) is
\beq{co2}
({\rm ad}_q \xi)\;\ov{T}^{\,I}_m =\sum\limits_n \Bigl(\rho^I(S(\xi))
\Bigr)_{mn}\,\ov{T}_n^{\,I}.
\eeq

\vspace*{2mm}

\red {\bf b) Quantum algebras in $R$-matrix formalism.}
\vspace*{2mm}

Let ${\cal J}_q=U_q({\cal J})$ be a quantum Lie algebra \cite{Dr,J,FRT} with
the unit element $e$, the co-product $\de$, the antipode $S$ and the co-unit
$\epsilon$. We suppose that the algebra ${\cal J}_q$ under consideration is
{\it quasitriangular} \cite{Dr}, i.e., there exists an invertible element
(called universal $R$-matrix) ${\cal R}\in{\cal J}_q \otimes {\cal J}_q$
with the properties:
\beq{qt1}
{\cal R}\,\de (\xi) \,=\, {\de}^{\prime} (\xi) \, {\cal R} \;\;
{\rm \ \ for \ any \ } \xi \in {\cal J}_q ,
\eeq
\beq{qt2}
(id\otimes\de)\,{\cal R} = {\cal R}_{13} {\cal R}_{12}, \;\;\;
(\de\otimes id)\,{\cal R} = {\cal R}_{13} {\cal R}_{23},
\eeq
where ${\de}^{\prime}=\sigma\circ{\de}$ and $\sigma$ denotes a permutation
in ${\cal J}_q \otimes {\cal J}_q\,$:
$\sigma(\xi\otimes\ti{\xi})=\ti{\xi}\otimes \xi$.
{}From (\ref{qt2}) one can derive the Yang-Baxter equation
\beq{yb}
{\cal R}_{12} {\cal R}_{13} {\cal R}_{23}= {\cal R}_{23} {\cal R}_{13}
{\cal R}_{12}.
\eeq
An action of the antipode $S$ on the universal $R$-matrix is given by
\beq{ant}
(S\otimes id) \, {\cal R}= (id\otimes S^{-1})\, {\cal R}= {\cal R}^{-1}.
\eeq

Let $\rho^I$ be a certain representation of ${\cal J}_q$. Then the objects
called $L$-operators are defined as follows
\beq{lop}
L_+^I = (\rho^I \otimes id) \, {\cal R}, \;\;\;
L_-^I = (\rho^I \otimes id) \, ({\cal R}^{\prime})^{-1},
\eeq
where ${\cal R}^{\prime}=\sigma\,(\,{\cal R}\,)$. Since
$L_{\pm}^I\in {\rm End}\,V^I\!\otimes\!{\cal J}_q$, they can be regarded
as matrices. From (\ref{yb}) the Yang-Baxter equations for the $L$-operators
are derived
\beq{yb2}
\R{}{IJ}{\pm}\,\Lo{1}{I}{\pm}\,\Lo{2}{J}{\pm} = \Lo{2}{J}{\pm}\,
\Lo{1}{I}{\pm}\, \R{}{IJ}{\pm},\;\;\;\;\:
\R{}{IJ}{+}\,\Lo{1}{I}{+}\,\Lo{2}{J}{-}
 = \Lo{2}{J}{-}\,\Lo{1}{I}{+}\,\R{}{IJ}{+},
\eeq
where $\Lo{1}{I}{\pm}\equiv {L^I}_{\pm} \!\otimes\! e^J
\in {\rm End}\,V^I\!\otimes\! {\rm End}\,V^J\!\otimes\!{\cal J}_q$ etc.
and $R_{\pm}$ are defined as follows
\beq{R}
\R{}{IJ}{+}=(\rho^I\otimes\rho^J)\,{\cal R},\;\;\;
\R{}{IJ}{-}=(\rho^I\otimes\rho^J)\,({\cal R}^{\prime})^{-1}.
\eeq
Further we shall omit indices of representations $I$, $J$ of $L$-operators
and $R$-matrices when it does not lead to uncertainty.

{\it Remark}. The matrix $L$ constructed from $L_+$, $L_-$ as follows
\beq{11}
  L = L_+ \, L_-^{-1},
\eeq
satisfies the following relation \cite{RS} (here $L^I=L_+^I (L_-^I)^{-1}$)
\beq{12}
 \Lo{1}{I}{} (R_-^{IJ})^{-1} \Lo{2}{J}{} R_-^{IJ} =
 (R_+^{IJ})^{-1} \Lo{2}{J}{} R_+^{IJ} \Lo{1}{I}{},
\eeq
which is known as "reflection equation" in the theory of
factorizable scattering \cite{re}.

Applying the representation $\rho^I$ to (\ref{qt2}) and (\ref{ant}), one
defines the Hopf algebra operations for the $L$-operators (see \cite{FRT} for
details)
\beq{Hopf}
\de \, L_{\pm} \, = \, L_{\pm} \, \dot{\otimes} \, L_{\pm}, \;\;\;\;
S(L_{\pm}) \,=\, L_{\pm}^{-1},
\eeq
where $\dot{\otimes}$ means a matrix product in the auxiliary space
$V^I$ and a tensor product in the quantum space ${\cal J}_q$. Relations
(\ref{Hopf}) rewritten for entries of the $L$-operators look like
\beq{Hopf2}
\de \,(L_{\pm})_{ij}\,=\, \sum\limits_{k} (L_{\pm})_{ik}\, {\otimes}\,
(L_{\pm})_{kj},\;\;\;\; S((L_{\pm})_{ij})\,=\,(L_{\pm}^{-1})_{ij}.
\eeq
\vspace*{2mm}

\red {\bf c) The case of $U_q(sl(n))$.}
\vspace*{2mm}

To exemplify the definitions given above we consider the case of
${\cal J}_q=U_q(sl(n))$ which will be our main example in the sequel.
The universal $R$-matrix for $U_q(sl(n))$ is given by \cite{BR}
\beq{univ}
{\cal R} = q^{2 \sum\limits_{ij} a_{ij}^{-1}\,H_i\otimes H_j}\,
{\prod}^{<}_{\al\in\Phi^{+}} \,E_{q^{-2}}\left( (1-q^{-2})\,
(q^{H_{\al}}X^+_{\al})\otimes (q^{-H_{\al}}X^-_{\al}) \right),
\eeq
where
$$ E_q(x)=\sum\limits_{r=0}^{\infty}\frac{x^r}{[r;q]!},
\;\;\;\;[r;q]!=\prod_{k=1}^r \,\frac{(q^k-1)}{(q-1)}$$
is a $q$-exponential, $a_{ij}$ is the Cartan matrix, $\Phi^+$ is the set of
positive roots and the factors in (\ref{univ}) are ordered in a specific way.

The commutation relations and the Hopf algebra attributes for all generators
of $U_q(sl(n))$ can be extracted from the matrix equations (\ref{yb2}) and
(\ref{Hopf}), where $L_{\pm}$ and $R_{\pm}$ are constructed by formulae
(\ref{lop}), (\ref{R}) from the $R$-matrix (\ref{univ}) (see e.g.
\cite{FRT,BR}). In particular, for the generators corresponding to simple
roots one gets
\beq{Hopf3}
\ar{cc} \de \,q^{H_i} = q^{H_i} \otimes q^{H_i} , \; &
\de \,X_i^{\pm} = X_i^{\pm} \otimes q^{-H_i} + q^{H_i} \otimes X_i^{\pm}, \;
\\ & \\ S(H_i) = -H_i, \; & S(X_i^{\pm}) = -q^{\mp 1} X_i^{\pm},
 \; \er \;\;\; i=1,..,n-1.
\eeq

Substituting (\ref{Hopf3}) into the general definition (\ref{6}) and taking
into account that the matrices $\rho(H_i)$ are diagonal in the standard
basis, one gets the definition of (covariant) tensor operators for
$U_q(sl(n))$
\beq{6b}
\ar{c}
q^{H_i}\, T_m^{\,J}\, q^{-H_i} =
q^{\left(\rho^J(H_i)\right)_{mm}}\, T_m^{\,J}, \\  \\
X_i^{\pm} \, T^{\,J}_m q^{H_i} - q^{H_i\mp 1}\, T^{\,J}_m \, X_i^{\pm} \, =
\sum\limits_n \;
T^{\,J}_n \left(\rho^J (X_i^{\pm})\right)_{nm} , \er \;\;\; i=1,..,n-1.
\eeq

\vspace*{3mm}

\red {\large\bf 2. Tensor operators in $R$-matrix approach}
\vspace*{2mm}

\red {\bf a) Generating matrices.}
\vspace*{2mm}

Let ${\cal J}_q$ be a deformed quasitriangular Lie algebra with the
corresponding universal $R$-matrix $\cal R$, let ${\cal M}$ be the
corresponding model space and let $\rho^I$ be a certain finite dimensional
representation of ${\cal J}_q$. Let us introduce new objects
\beq{u}
  U^J,\;W^J\;\;\in\, {\rm End}\,V^J \otimes {\rm End}\,{\cal M},
\eeq
which can be regarded as matrices in the auxiliary space $V^J$
with entries being operators mapping the model space  ${\cal M}$ into itself.

\begin{prop}
Let ${\rho}^{I_0}$ and $\rho^J$ denote the fundamental and arbitrary
finite dimensional representations of the quantum Lie algebra ${\cal J}_q$,
respectively. Suppose that $U^J$ introduced in (\ref{u}) satisfies the
relations
\beq{2.1}
\Lo{1}{I_0}{+}\,\Uo{2}{J}=\,\Uo{2}{J}\,\R{}{I_0 J}{+}\,\Lo{1}{I_0}{+},
\;\;\;\;\; \Lo{1}{I_0}{-}\,\Uo{2}{J}\,=\,\Uo{2}{J}\,\R{}{I_0 J}{-} \,
\Lo{1}{I_0}{-} ,
\eeq
where $\Lo{}{I_0}{\pm}$, $\R{}{I_0 J}{\pm}$ are defined as in (\ref{lop}),
(\ref{R}). Then $n$-th row (for any $n$) of the matrix $U^J$ forms a
$q$-tensor operator (which transforms by the representation $\rho^J$) in
the sense of definition (\ref{6}) according to the rule:
\beq{2.10}
U^J_{ni}= T_{\;i}^{J\,\langle n\rangle},\;\;\; i=1,.., {\rm dim}\,\rho^J .
\eeq
\end{prop}

Below we call the matrix $U^J$ {\it generating matrix}~ because,
according to the Wigner-Eckart theorem, its entries produce
corresponding CGC.

{\it Remark}. Obviously, once we have the equations (\ref{2.1}) with the
fundamental $L$-operators $L_{\pm}^{I_0}$, we can extend them to the case of
$L_{\pm}^I$ corresponding to the arbitrary representation $\rho^I$ as well.
These equations will describe $q$-tensor operators in terms of polynomial
combinations of the generators of the deformed Lie algebra ${\cal J}_q$.

{\it Remark}. An equivalent form of the relations (\ref{2.1}) (due to the
previous Remark, we can write it down for arbitrary representation of the
$L$-operators) is
\beq{2.2}
\Lo{1}{I}{+}\, \Uo{2}{J}\, (\Lo{1}{I}{+})^{-1} \,=\,\Uo{2}{J}\,
\R{}{I J}{+}, \;\;\;\; \Lo{1}{I}{-}\,\Uo{2}{J}\,
(\Lo{1}{I}{-})^{-1}\,= \,\Uo{2}{J}\,\R{}{I J}{-} .
\eeq

{\it Remark}. The consequence of (\ref{2.1}) for $L$ defined in (\ref{11}) is
\[ \Lo{1}{I}{}\;\Uo{2}{J}\,\R{}{IJ}{-}=\Uo{2}{J}\,\R{}{IJ}{+}\,\Lo{1}{I}{}.\]

{\it Proof} of the Proposition 1.
Let us notice that the l.h.s. of (\ref{2.2}) is nothing but the $R$-matrix
form of the definition of adjoint action (\ref{5}).
Indeed, let $L_{\pm}$ be N$\times$N matrices (we omit the index $I_0$ since
the following arguments hold for any representation $\rho^I$). Then, using
the matrix form of tensor product
$$\up{1}A\;\up{2}B=\!\left(\!\ar{ccc}\scriptstyle A_{11}B &\!\! \dots \!\! &
\scriptstyle A_{1N}B \\ \scriptstyle \vdots & &\scriptstyle \vdots \\
\scriptstyle A_{N1}B & \!\! \dots \!\!& \scriptstyle A_{NN}B \er\!\!\right),$$
one can write down the product
$\Lo{1}{}{\pm}\,\Uo{2}{J}(\Lo{1}{}{\pm})^{-1}$
as N$\times$N matrix with entries (which are matrices as well) given by
\beq{a1}
 \Bigl(\up{1}{L}_{\pm}\,\Uo{2}{J} (\up{1}{L}_{\pm})^{-1}\Bigr)_{ij}\,=\,
\sum\limits_{k}\,\bigl(L_{\pm}\bigr)_{ik}\,U^J\,S\Bigl(\bigl(L_{\pm}
\bigr)_{kj}\Bigr)\,\equiv\, \left({\rm ad}_q\,\bigl(L_{\pm}\bigr)_{ij}\right)
\, U^J,\;\;\; i,j=1,..,N.
\eeq
Here we have used the relations (\ref{Hopf2}) and the definition (\ref{5}).
The r.h.s. of (\ref{2.1}) written in the same manner looks like
\beq{a2}
\Bigl(\Uo{2}{J} \R{}{IJ}{\pm} \Bigr)_{ij}\,=\,
U^J\;\Bigl( (id\otimes\rho^J)\,L_{\pm}^I \Bigr)_{ij},
\eeq
where we have used the identities $R_{\pm}^{IJ}=(id\otimes\rho^J)L_{\pm}^I$
which follow from (\ref{lop}) and (\ref{R}).

{}From (\ref{2.2})-(\ref{a2}) it follows that each row of the matrix $U^J$
obeys the vector form (\ref{2.4b}) of the definition (\ref{6}). Since the
pair $L_{\pm}^{I_0}$ comprises all the generators of ${\cal J}_q$, we
conclude that each row of $U^J$ is a $q$-tensor operator. \qed

{\it Remark}. In the case of arbitrary quantum algebra one should have the
equations (\ref{2.1}) for the sufficient set of $L$-operators (not only for
the fundamental ones).

{\it Remark}.
Multiplying (\ref{2.1}) from the left by $e\otimes\vec{x}^{\,t}_r$
($e$ stands for the unit element of ${\cal J}_q$), where
$(\vec{x}_r)_k=\dl_{rk}$, $k\!=\! 1,..,N$, one gets equations of the same
type for $r$-th row of $U^J$. This implies that the rows of $U^J$ are
independent, i.e., in the Proposition 1 one may regard $U^J$ not as a
{\it matrix} but as a {\it vector} (row) defined as $U^J\in V^J\otimes
{\rm End}\,{\cal M}$. But, bearing in mind possible applications, it is
preferable to arrange these vectors in a square matrix. The most interesting
case is when all the rows correspond to nonequivalent \footnote{ Tensor
operators are equivalent if their matrix elements are proportional.}
tensor operators.

{\it Remark}. From the previous Remark it also follows that an arbitrary
linear combination of the rows of $U^J$ defines some $q$-tensor operator as
well.

{\it Remark}. One should not be confused by the fact that equations (\ref{2.1})
with $U^J$ replaced by $g\in{\rm Fun}\,(G_q)=\left(U_q({\cal J})\right)^*$
coincide with the well known duality relations between the quantum algebra and
the quantum group. These equations deal with essentially different objects.

Now let us define conjugated (contravariant) tensor operators.
\begin{prop}
Let $\Lo{}{I_0}{\pm}$, $R_{\pm}^{I_0 I}$, $W^J$ be defined as in (\ref{lop}),
(\ref{R}), (\ref{u}). Suppose that $W^J$ satisfies the relations
\beq{co3}
\Lo{1}{I_0}{+}\,\W{2}{J}\,=(R_+^{I_0 J})^{-1}\,\W{2}{J}\,
\Lo{1}{I_0}{+},\;\;\;\;\; \Lo{1}{I_0}{-}\,\W{2}{J}=
(R_-^{I_0 J})^{-1}\,\W{2}{J}\,\Lo{1}{I_0}{-}.
\eeq
Then $n$-th column (for any $n$) of the matrix $W^J$ forms a conjugated
$q$-tensor operator (which transforms by the representation $\rho^J$) in
the sense of definition (\ref{co2}) according to the rule:
\beq{co4}
 W^J_{in} = \ov{T}_{\;i}^{J\,\langle n\rangle},\;\;\; i=1,..,
{\rm dim}\,\rho^J .
\eeq
\end{prop}

{\it Proof}$\,$ is analogous to the one that was given for the Proposition 1.
One should use the identities $(R_{\pm}^{I_0 J})^{-1}=(id\otimes\rho^J)\,
S(L_{\pm}^J)$ which follow from (\ref{ant}),(\ref{lop}) and (\ref{R}). \qed

{\it Remark}. An equivalent form of the relations (\ref{co3}) is
\beq{co5}
\Lo{1}{I}{+}\,\W{2}{J}(\Lo{1}{I}{+})^{-1}\,=(\R{}{IJ}{+})^{-1}\,\W{2}{J},
\;\;\;\; \Lo{1}{I}{-}\,\W{2}{J} (\Lo{1}{I}{-})^{-1}\,=(\R{}{IJ}{-})^{-1}\,
\W{2}{J}.
\eeq

{\it Remark}. The consequence of (\ref{co3}) for $L$ defined in (\ref{11}) is
\[ \R{}{IJ}{+} \Lo{1}{I}{+}\,(\R{}{IJ}{-})^{-1} \W{2}{J} =
\W{2}{J}\, \Lo{1}{I}{+}. \]

{\it Remark}. It is easy to see that if the matrix $U^J$ corresponds to
$q$-tensor operator (i.e., obeys (\ref{2.1})) and is invertible, then
$W^J=(U^J)^{-1}$ corresponds to conjugated $q$-tensor
operator (i.e., it obeys (\ref{co3})) and vice versa.
\vspace*{2mm}

Now let us remind that for classical groups the "scalar product"
$\sum\limits_m\,Q_m\,\ov{T}_m$ of covariant and contravariant tensor operators
is a scalar with respect to the adjoint action of the algebra \cite{0}.
The $R$-matrix approach allows to formulate and prove this statement for a
quantum algebra ${\cal J}_q$ as follows.
\begin{prop}
Let $U^J$ and $W^J$ be generating matrices obeying (\ref{2.1})
and (\ref{co3}), respectively. Then all the entries of the matrix
$Z^J=U^J\,W^J$ are scalars with respect to the adjoint action of
${\cal J}_q$.
\end{prop}
{\it Proof}. First, from (\ref{2.1}) and (\ref{co3}) one finds
\beq{sc}
\ar{c} \Lo{1}{I_0}{\pm}\,\Z{2}{J}\,(\Lo{1}{I_0}{\pm})^{-1}=
(\Lo{1}{I_0}{\pm}\,\Uo{2}{J}\,(\Lo{1}{I_0}{\pm})^{-1})\;
(\Lo{1}{I_0}{\pm}\,\W{2}{J}\,(\Lo{1}{I_0}{\pm})^{-1})= \\ \\
= \, (\Uo{2}{J} \R{}{I_0 J}{\pm})\,((\R{}{I_0 J}{\pm})^{-1} \W{2}{J})
= \, \Uo{2}{J}\,\W{2}{J} = \,\Z{2}{J}.  \er
\eeq
Next, multiplying (\ref{sc}) by $e\otimes\vec{x}^{\,t}_l$ from the left
and by $e\otimes{\vec{x}}_r$ from the right (these number-valued vectors
$(\vec{x}_r)_k=\dl_{rk}$ commute with $\up{1}{L}$), we get
this relation for each element $Z^J_{\,lr}$. Finally, using the same
arguments as in the proof of the Proposition 1, we conclude that each
$Z^J_{\,lr}$ satisfies (\ref{6}) with a trivial (one-dimensional)
representation $\rho^0$ (which is nothing but the co-unit $\epsilon$) in the
r.h.s, and, therefore, it is a scalar. \qed

{\it Remark}. In other words, the equation (\ref{sc}) for each $Z^J_{\,lr}$
may be regarded as (\ref{2.2}) or (\ref{co5}) with the trivial $R$-matrix
$R_{\,\pm}=e^{I_0}\otimes 1$.
\vspace*{2mm}

\red {\bf b) Wigner-Eckart theorem.}
\vspace*{2mm}

Let us give a proof of the Wigner-Eckart theorem for $q$-tensor operators with
the help of the $R$-matrix language. We restrict our consideration to the case
of quantum simple quasitriangular Lie algebra ${\cal J}_q$. We shall prove the
following equivalent of the statement (\ref{WE})
\beq{t1}
  T^J_m\,|Im^{\prime}\rangle\,=\,\sum\limits_K ||T^J||_{I\,K}\,
\Bigl\{\!\ar{ccc} I &\!J\!&\!K\\ m^{\prime}\!&\! m\!&\!
m^{\prime\prime}\er\!\Bigr\}\,|K m^{\prime\prime}\rangle.
\eeq

In this section we shall use the equation (\ref{2.1}) with $U^J$ being not
a matrix but a single row (as explained above, all the rows are independent).
So, let $U^J$ be a row which forms the tensor operator $T^J$. Then (\ref{t1})
can be rewritten as follows ( $t$ denotes transposition)
\beq{t2}
 \Uo{2}{J}\,\ve{I\,}{1}{\;t}\,=\,\sum\limits_K \, ||T^J||_{I\,K}\,
 \ve{K}{}{\!t} \,{\rm C}[IJK],
\eeq
where $\ve{K}{}{}$ denotes the vector of basic vectors for given subspace
${\cal H}^K$, i.e., \hbox{ $|K\rangle_m\!=\!|Km\rangle$ } (compare with
(\ref{2.4a})) and we introduced the Clebsch-Gordan maps \cite{KR,MS}
defined as
\beq{t3}
{\rm C}[IJK] :\, V^I\otimes V^J \rar V^K,\;\;\;\;
{\rm C}'[IJK] :\,V^K\rar V^I\otimes V^J .
\eeq
These objects have the properties
\beq{t4}
\ar{c}  {\rm C}[IJK]\, (\rho^I\otimes \rho^J)\,\de(\,\xi\,)=
\rho^K(\xi) \, {\rm C}[IJK],  \\ [2mm]
 (\rho^I\otimes \rho^J)\,\de(\,\xi\,)\,{\rm C}'[IJK]=
 {\rm C}'[IJK]\,\rho^K(\xi) , \er \ {\rm for\ any\ } \xi\in{\cal J}_q ,
\eeq
\beq{t5}
 \sum\limits_K {\rm C}'[IJK] \,{\rm C}[I^{\prime} J^{\prime} K]
= \dl_{II^{\prime}}\,\dl_{JJ^{\prime}},\;\;\;\;\;
{\rm C}[IJK]\,{\rm C}'[IJL]=\dl_{KL}
\eeq
and for fixed representations they can be regarded as rectangular matrices
(with entries being numerical values of the CGC),
e.g. the r.h.s. of (\ref{t5}) are the identity matrices.

{}From (\ref{t4})-(\ref{t5}) and the quasitriangularity relation (\ref{qt2})
one derives the formulae which are known as the {\it fusion procedure} for
$R$-matrices:
\beq{t6}
 \up{23}{\rm C}[IJK] \, \R{13}{LJ}{}\,\R{12}{LI}{}\,
=\,\R{}{LK}{}\;\up{23}{\rm C}[IJK], \;\;\;\;
 \R{13}{LJ}{}\,\R{12}{LI}{}\,\up{23}{\rm C'}[IJK] \,
=\,\up{23}{\rm C'}[IJK] \R{}{LK}{},
\eeq
where we denoted $\up{23}{\rm C}[IJK]\equiv id\otimes {\rm C}[IJK]$,
$\up{23}{\rm C'}[IJK]\equiv id\otimes {\rm C'}[IJK]$.

{\it Remark}. The relations (\ref{t6}) look very similarly to (\ref{2.1}).
But we have to stress that these relations describe essentially different
objects.

Next, we notice that, due to (\ref{lop}), action of generators of
${\cal J}_q$ on basic vectors can be written as follows
\beq{t7}
   \Lo{1}{I}{\pm}\,\ve{K}{2}{\!t} \,=\,\ve{K}{2}{\!t}\,\R{}{IK}{\pm}.
\eeq

Now we are able to prove the Wigner-Eckart theorem. Indeed, using (\ref{2.1})
and (\ref{t6})-(\ref{t7}), we find
\[ \ar{l} \Lo{1}{N}{\pm}\,\Bigl(\Uo{3}{J}\,\ve{I\,}{2}{\;t}\,
\up{23}{{\rm C'}}[IJK]\Bigr)=
\Uo{3}{J}\, \R{13}{NJ}{\pm}\, \Lo{1}{N}{\pm}\,
\ve{I\,}{2}{\;t}\, \up{23}{{\rm C'}}[IJK]\,= \\  \\ =
\Uo{3}{J}\, \R{13}{NJ}{\pm}\,\ve{I\,}{2}{\;t}\,
\R{12}{NI}{\pm}\,\up{23}{{\rm C'}}[IJK]\,=
\,\Uo{3}{J}\,\ve{I\,}{2}{\;t}\, \R{13}{NJ}{\pm}\,
\R{12}{NI}{\pm}\,\up{23}{{\rm C'}}[IJK]\,=  \\ \\
= \Bigl( \Uo{3}{J}\,\ve{I\,}{2}{\;t}\,\up{23}{{\rm C'}}[IJK]\Bigr)\,
\R{}{NK}{\pm}. \er \]
According to (\ref{t7}), this implies that the vector
$(\Uo{2}{J}\ve{I\,}{1}{\;t} {\rm C'}[IJK])$ coincides (up to an arbitrary
constant depending on $I$, $J$, $K$) with the vector $\ve{K}{}{\!t}$ :
\[ \Uo{2}{J}\,\ve{I\,}{1}{\;t}\, {\rm C'}[IJK]\,=\,
||T^J||_{IK}\, \ve{K}{}{\!t} .\]
Due to (\ref{t5}) the later equality is equivalent to (\ref{t2}).
Thus, the Wigner-Eckart theorem is proved.

{\it Remark}. If the algebra ${\cal J}$ is not simple, the theory of tensor
operators is more sophisticated because multiplicities appear in tensor
products of representations. In particular, it makes a proof of the
Wigner-Eckart theorem more complicated.
One may assume that the information about these multiplicities is hidden
in the structure of generating matrices and $R$-matrices.
\vspace*{2mm}

\red {\bf c) Fusion procedure}
\vspace*{2mm}

The relations (\ref{t6}) describe the fusion procedure for $R$-matrices,
i.e., they allow in principle to construct from the fundamental $R$-matrix
$R^{I_0 I_0}$ all other representations $R^{IJ}$ of the universal $R$-matrix
for a given quantum Lie algebra. In what follows we shall use slightly
different form of fusion formula (originally introduced in \cite{KS} for
spectral-dependent $R$-matrices)
\beq{f1}
  \R{1,32}{\,LK}{\,\pm} =\Po{23}{IJ}{K}\,\R{13}{LJ}{\pm}\,\R{12}{LI}{\pm}\,
 \Po{23}{IJ}{K},
\eeq
where the l.h.s. stands for $R^{\,LK}_{\pm}$
written in the basis of the bigger space \hbox{$V^L\otimes V^I\otimes V^J$}
and $P^{IJ}_{K}$ denotes a projector (i.e., $(P^{IJ}_{K})^2=P^{IJ}_{K}$)
onto the subspace in $V^I\otimes V^J$ corresponding to the representation
$\rho^K$.

{\it Remark.} For example, such projectors appear in the decomposition
formula for the fundamental $U_q(sl(n))$ $R$-matrix (see \cite{FRT} for
other principal series)
\beq{f2}
  \wh{R}_{\pm}=q^{\mp \frac{1}{n} }\,\left( q^{\pm 1}\,P_+
- q^{\mp 1}\,P_- \right),
\eeq
where $(\wh{R}_{\pm})_{ij,kl}=({R}_{\pm})_{ji,kl}$ and $P_{\pm}$ are the
projectors in ${\bf C}^n\otimes{\bf C}^n$ ($q$-symmetrizer and
$q$-antisymmetrizer) of ranks $\frac{n(n+1)}{2}$ and $\frac{n(n-1)}{2}$,
respectively.

The statement that (\ref{f1}) satisfies the Yang-Baxter relation is equivalent
\cite{KS} to the following property of the projectors
\beq{f21}
\Po{23}{IJ}{K} \,\R{13}{LJ}{\pm}\,\R{12}{LI}{\pm}\,=\,\R{13}{LJ}{\pm}\,
 \R{12}{LI}{\pm}\,\Po{23}{IJ}{K}.
\eeq
Notice that in the case of $U_q(sl(n))$ for $\rho^I\!=\!\rho^J\!=\!\rho^{I_0}$
this property obviously follows from the decomposition (\ref{f2}). The general
case is discussed in the Appendix.

{\it Remark}. In the case of $U_q(sl(2))$ one can construct \cite{KS} all the
matrices $R^{IJ}_{\pm}$ from the fundamental ones
$R^{{\frac 12}\,{\frac 12}}_{\pm}$ with the help only of the symmetrizers
$P_+^s$ which are the projectors in $({\bf C}^2)^{\otimes 2s}$ onto
$(2s+1)$-dimensional subspace.

Now we are going to develop the fusion procedure for generating matrices.
\begin{prop}
Let $U_I$, $U_J$ and $W_I$, $W_J$ satisfy (\ref{2.1})
and (\ref{co3}), respectively, and let $\F{}{IJ}\in V^I\otimes V^{J}$ be a
matrix with entries commuting with all the generators of ${\cal J}_q$. Then
\beq{f3}
\Uf{21}{IJ}{K} =\F{12}{IJ}\,\Uo{2}{J}\,
\Uo{1}{I}\,\Po{12}{IJ}{K},\;\;\;\; {\rm and} \;\;\;\;
\Wf{12}{\,IJ}{K} = \Po{12}{IJ}{K} \,\W{1}{I}\,
\W{2}{J}\, \F{12}{IJ}
\eeq
obey (\ref{2.1}) and (\ref{co3}), respectively, with $R$-matrices
constructed by the formulae (\ref{f1}).
\end{prop}

{\it Proof}. Using (\ref{2.1}) and (\ref{f1}), we verify this statement for
$U^{IJ}_K$ (for $W^{IJ}_K$ it can be done in
the same way)
\[ \ar{c} \Lo{1}{N}{\pm} \,\Uf{32}{IJ}{K}
=\,\Lo{1}{N}{\pm}\, \F{23}{IJ} \,\Uo{3}{J}\,
\Uo{2}{I} \,\Po{23}{IJ}{K} =\\ \\
= \F{23}{IJ} \, \Uo{3}{J} \, \R{13}{NJ}{\pm}
\,\Lo{1}{N}{\pm}\, \Uo{2}{I} \, \Po{23}{IJ}{K} =
\F{23}{IJ} \,\Uo{3}{J} \,\R{13}{NJ}{\pm}\,\Uo{2}{I} \,
\R{12}{NI}{\pm}\,\Lo{1}{N}{\pm}\,\Po{23}{IJ}{K}= \\ \\
= \F{23}{IJ} \, \Uo{3}{J} \, \Uo{2}{I} \,
\Po{23}{IJ}{K} \, \R{13}{NJ}{\pm} \,\R{12}{NI}{\pm}\,
\Po{23}{IJ}{K} \Lo{1}{N}{\pm} \,=\,
\Uf{32}{IJ}{K} \,\R{1,32}{NK}{\pm} \Lo{1}{N}{\pm} .\er \]
Here we have used the property (\ref{f21}) and formula (\ref{f1}) and
exploited the commutativity of all entries of $F$ with those of $L_+$,
$L_-$. \qed

Starting with the fundamental generating matrices and applying the
described procedure appropriately many times, one can construct
generating matrices corresponding to the arbitrary representation.

{\it Remark}. Due to the presence of the matrix $F$ in (\ref{f3}), these
formulae produce families of tensor operators $U^{JJ'}_K$ for given
$J$, $J'$, $K$. However, one will find only finite number of nonequivalent
tensor operators among them.

{\it Remark}. It is important for possible applications (see the Discussion
below) that the matrix $F$ may depend nontrivially on the central
elements of the algebra ${\cal J}_q$.

{\it Remark}. Formulae (\ref{f3}) give the generating matrices $U^K$, $W^K$
in the basis of \hbox{ $V^I\otimes V^J$ }. One can rewrite them in the
"natural" basis of the space $V^K$ as follows
\[ U^K_{\,ij}\,=\,e_i^*\; U^{JJ^{\prime}}_K \; e_j , \;\;\;\;
i,j=1,..,N\equiv {\rm dim}\,\rho^K. \]
Here $e_i$ are the eigenvectors of the projector $P^{JJ^{\prime}}_K$ and
$e^*_i\,e_j=\dl_{ij}$.

The fusion procedure described above allows also to construct the invariant
operators of ${\cal J}_q$. For simplicity, we consider only the case of
$U_q(sl(n))$.
\begin{prop}
Let $U^{I_0}$ and $W^{I_0}$ obey, respectively, (\ref{2.1}) and (\ref{co3})
for the fundamental representation $\rho^{I_0}$ of ${\cal J}_q=U_q(sl(n))$,
let $P^{I_0}_{-}$ denote the antisymmetrizer in $(V^{I_0})^{\otimes n}$ and
let $\{ F_k\}_{k=1}^{n-1}$ be a set of matrices such that
$F_k=f_k\otimes e^{I_0}...\otimes e^{I_0}\in (V^{I_0})^{\otimes n}$,
where entries of the matrices $f_k\in (V^{I_0})^{\otimes (k+1)}$ commute with
all the elements of ${\cal J}_q$. Then all the entries of the matrices
\[ \Uf{n..1}{0}{}=F_{n-1}\,\Uo{n}{I_0}\,...\,F_1\,\Uo{2}{I_0}\,
\Uo{1}{I_0}\, P^{I_0}_{-} \ \ {\rm and} \ \
\Wf{1..n}{\;0}{}=P^{I_0}_{-}\, \W{1}{I_0} \W{2}{I_0} F_1\,...
\W{n}{I_0} F_{n-1}  \]
commute with any $\xi \in {\cal J}_q$. 
\end{prop}

{\it Proof}. For $\Uf{}{0}{}$ one finds (the case of $W^0$ is
analogous)
\bea{}
 \Lo{0}{I_0}{\pm} \, \Uf{n..1}{0}{} &=&
\Lo{0}{I_0}{\pm} \, F_{n-1}\,\Uo{n}{I_0}\,...\,F_1\,\Uo{2}{I_0}\,
\Uo{1}{I_0}\,P^{I_0}_{-}\, =  \nn \\
&=&\, F_{n-1}\,\Uo{n}{I_0}\,\R{n0}{I_0I_0}{\pm}\,...\,F_1\,\Uo{2}{I_0}\,
\R{20}{I_0I_0}{\pm}\,\Uo{1}{I_0}\R{10}{I_0I_0}{\pm}\,\,P^{I_0}_{-}\,
\Lo{0}{I_0}{\pm} = \nn \\
&=& F_{n-1}\,\Uo{n}{I_0}\,...\,F_1\,\Uo{2}{I_0}\,\Uo{1}{I_0}\,
\R{n0}{I_0I_0}{\pm}\,...\,\R{10}{I_0I_0}{\pm}\, P^{I_0}_{-}\,
\Lo{0}{I_0}{\pm}\,= \nn \\
&=& F_{n-1}\,\Uo{n}{I_0}\,...\,F_1\,\Uo{2}{I_0}\,\Uo{1}{I_0}\,
P^{I_0}_{-}\,\Lo{0}{I_0}{\pm}\,=\, \Uf{n..1}{0}{}\,\Lo{0}{I_0}{\pm}. \nn
\eea
Here we have used the property \cite{KS,J,FRT}:
$\R{n0}{I_0I_0}{\pm}\,...\,\R{10}{I_0I_0}{\pm}\,P^{I_0}_{-}=c\,P^{I_0}_{-}$
(we have $c=1$ due to the special choice of normalization in (\ref{f2})).
\qed

{\it Remark}. The construction described in the Proposition 5 is a
modification of the formula for the {\it quantum determinant} \cite{KS,FRT}.
Hence it can be generalized to the case of other quantum Lie algebras as well.
\vspace*{2mm}

\red {\bf d) Connection of covariant and contravariant operators.}
\vspace*{2mm}

It is well known \cite{0} that entries of contravariant tensor operators
can be reexpressed via those of covariant tensor operators and vice versa.
Let us formulate this correspondence in the $R$-matrix language. Recall that
for $U_q(sl(n))$ there exists the $q$-Weyl element ${\cal W}$ \cite{KR,RT}
(in general, it belongs to the extension of $U_q(sl(n))$) with the property
\beq{w1}
{\cal W}\, \xi\, {\cal W}^{-1} \, = \, \tau\,\left( S(\xi)\right)
\;\;\;\;{\rm for \ any \ } \xi \in {\cal J}_q,
\eeq
where $\tau$ is an antiautomorphism defined as follows
\beq{w3}
  \tau(X_{i}^{\pm}) = X_{\theta(i)}^{\mp},\;\;\; \tau(H_{i})=H_{\theta(i)},
\;\;\;\; i=1,..,n-1.
\eeq
Here $\theta$ stands for the certain involution of Dynkin diagram (it can be
extracted from (\ref{univ})) such that $\theta\!\circ\!\theta\!=\!id$. It can
be represented as an element of the group of permutations (of positive roots).
For example: for $sl(2)$, $\theta\!\sim\!\left\{1\atop 1\right\}$;\ \ for
$sl(3)$, $\theta\!\sim\!\left\{{123}\atop {213}\right\}$.

Taking into account the relation $\tau S\tau=S^{-1}$, one gets from (\ref{w1})
two "crossing-symmetry" relations for the universal $R$-matrix
\beq{w2}
(\tau\otimes id)\,{\cal R}^{-1} = \up{1}{\cal W}\,{\cal R}\,
(\up{1}{\cal W})^{-1} ,\;\;\;\;\; (\up{2}{\cal W})^{-1}\,{\cal R}^{-1}
\up{2}{\cal W} =(id\otimes\tau){\cal R}.
\eeq
For each irreducible representation $\rho^I$ of $U_q(sl(n))$ the element
$\cal W$ can be written as a certain matrix ${\cal W}^I$. In the case of
$U_q(sl(2))$ the corresponding matrix for the irreducible representation of
spin $s$ is given by the (2s+1)$\times$(2s+1) matrix \cite{KR,RT}
\[ {\cal W}^{\,s}_{mk} = (-1)^{k}\,q^{s+1-m}\,\dl_{m,2s+2-k},
\;\;\;\; m,k=1,..,2s+1. \]
Let us give another example. Introducing n$\times$n matrices
\beq{E}
 \left(E_{i,j}\right)_{mn}\,=\,\dl_{im}\;\dl_{jn},
\eeq
one can write down the fundamental representation of $U_q(sl(n))$
as follows
\beq{e}
\rho^{I_0}(H_i)=E_{i,i}-E_{i+1,i+1},\;\;\;
\rho^{I_0}(X_i^{+})=E_{i,i+1}, \;\;\;
\rho^{I_0}(X_i^{-})=E_{i+1,i}, \;\;\; i=1,..,n-1.
\eeq
\begin{prop}
For the fundamental representation $\rho^{I_0}$ of $U_q(sl(n))$ the element
$\cal W$ is given by the following n$\times$n matrix
\beq{W}
{\cal W}^{I_0}_{mk}=(-1)^k \,q^{k-(n+1)/2}\,\dl_{m,n-k+1},
\;\;\;\;m,k=1,..,n.
\eeq
\end{prop}
{\it Proof}. From (\ref{E}), (\ref{W}) one derives
\beq{wew}
{\cal W}^{I_0}\,E_{i,j}\,\left(\,{\cal W}^{I_0}
\,\right)^{-1}=q^{i-j}\,(-1)^{i+j}\,E_{n-i+1,n-j+1} .
\eeq
Substituting (\ref{wew}) into (\ref{e}), one gets exactly the transformation
(\ref{w1}) for the generators corresponding to simple roots. It is
sufficient since (\ref{w1}) defines a homomorphism. \qed  \\
(Notice also that (\ref{wew}) gives an explicit form of the involution
$\theta$ in this case.)

Now the correspondence between covariant and contravariant tensor
operators can be formulated and proved in the $R$-matrix language as follows.
\begin{prop}
Let $U^J$ and $W^J$ be generating matrices obeying (\ref{2.1}) and (\ref{co3}),
respectively, and let ${\cal W}^J$ be a matrix form of the $q$-Weyl element
(\ref{w1}) in the representation $\rho^J$. Then
\[ \ti{U}^J\,\equiv\,(W^J)^t\,{\cal W}^J,\;\;\;\;\;
\ti{W}^J\,\equiv\,{\cal W}^J\, (U^J)^{t} \]
obey (\ref{2.1}) and (\ref{co3}), respectively ($t$ denotes the transposition
in auxiliary space).
\end{prop}
{\it Proof.} Let us prove, say, the first of these formulae. Taking
transposition of the equations (\ref{co3}) in the second auxiliary space,
multiplying them from the right by $(e\otimes{\cal W}^J)$ and applying
(\ref{w2}), one gets exactly the equations (\ref{2.1}) for $\ti{U}^J$. \qed

{\it Remark}. From the Propositions 3 and 7 we conclude that
$U^J{\cal W}^J (U^J)^t$ and \\
\hbox{ $(W^J)^t{\cal W}^J W^J$ } are the matrices of scalars.

Now let us consider the automorphism of $U_q(sl(n))$ which is given by
(note that $\tau$ in (\ref{w3}) was defined as an antiautomorphism)
\beq{ch1}
 \chi\,(X^{\pm}_i)\,=\,X^{\mp}_{\theta (i)},\;\;\;\;
 \chi\,({H}_i)\,=\,H_{\theta (i)}, \;\;\;\; i=1,..,n-1,
\eeq
where $\theta$ was described above. For any irreducible representation
$\rho^I$ it can be realized as a similarity transformation with the
antidiagonal matrix $\chi^I$
\beq{ch2}
 \rho^I(\chi\,(\xi))\,=\,\chi^I\,\rho^I(\xi)\,(\chi^I)^{-1},\;\;\;\;
\chi_{mk}^I=\dl_{m,N+1-k},\;\;\; m,k=1,..,N={\rm dim}\,\rho^I.
\eeq
It is easy to see that the fundamental $R$-matrices for $U_q(sl(n))$
\cite{Dr,J}
\beq{fr}
 R_{\,\pm}^{I_0 I_0} = q^{\pm \frac{1}{n}} \! \Bigl(q^{\pm 1}\,\sum\limits_i
E_{i,i}\otimes  E_{i,i} + \sum\limits_{i\neq j} E_{i,i}\otimes E_{j,j} \pm
\omega \!\!\! \sum\limits_{\pm (j-i) > 0} \!\! E_{i,j}\otimes
E_{j,i}\Bigr), \;\;\; i,j=1,..,n
\eeq
where $E_{i,j}$ were defined in (\ref{E}), have the properties
\beq{cr}
\ho{1}{I_0}\,R_{\pm}^{I_0 I_0}\,(\ho{1}{I_0})^{-1}=(R_{\pm}^{I_0 I_0})^{t_1},
\;\;\;\; \ho{2}{I_0}\,R_{\pm}^{I_0 I_0}\,(\ho{2}{I_0})^{-1}=
(R_{\pm}^{I_0 I_0})^{t_2}.
\eeq
\begin{prop}
Let $\rho^{I_0}$ denote the fundamental representation of $U_q(sl(n))$ and
let $\wh{U}^{I_0}$ and $\wh{W}^{I_0}$ be defined as in (\ref{u}) and satisfy
the relations
\beq{wu}
\Lo{1}{I_0}{\pm}\,\Uw{2}{I_0}\,=\,\R{}{I_0 I_0}{\pm}\,
\Uw{2}{I_0}\,\Lo{1}{I_0}{\pm}, \;\;\;\; \Lo{1}{I_0}{\pm}\,
\Ww{2}{I_0}\,=\,\Ww{2}{I_0}\, (\R{}{I_0 I_0}{\pm})^{-1}\,
\Lo{1}{I_0}{\pm},
\eeq
respectively. Then $n$-th column of $\wh{U}^{I_0}$ and $m$-th row (for any
$n$ and $m$) of $\wh{W}^{I_0}$ form co- and contravariant $q$-tensor operators
(which transform by the representation $\rho^{I_0}$), respectively,
according to the rules
\[ \wh{U}^{I_0}_{in}= T_{N-i}^{I_0\,\langle n\rangle},\;\;\;\;
 \wh{W}^{I_0}_{mi}= \ov{T}_{N-i}^{I_0\,\langle m\rangle},\;\;\;\;
i=1,..,N={\dim}\,\rho^{I_0}.\]
\end{prop}
{\it Proof}. Taking transposition of equations (\ref{2.1}) and (\ref{co3}) in
the second auxiliary space, multiplying them from the left and from the right,
respectively, by $(e\otimes\chi^{I_0})$ and using the identities (\ref{cr}),
one gets exactly the equations (\ref{wu}) for
$\wh{U}^{I_0}\equiv\chi^{I_0}\,(U^{I_0})^{t}$ and
$\wh{W}^{I_0}\equiv (W^{I_0})^{t}\,\chi^{I_0}$. \qed
\vspace*{2mm}

\red {\bf e) Classical limit.}
\vspace*{2mm}

Let us consider the limit $q\rar 1$ which returns us to the classical
representation theory. Let $r\in{\cal J}\otimes{\cal J}$ denote the classical
$r$-matrix, i.e., $r=\lim\limits_{q\rar 1}({\cal R}-e\otimes e)/(q-1)$.
For ${\cal J}=sl(n)$ it is given by \cite{Dr}
\beq{clr}
 r= 2\,(\sum\limits_{\al\in{\Phi}^+} H_{\al}\otimes H_{\al} +
X_{\al}^+\otimes X_{\al}^- ) .
\eeq
Similarly to the quantum case, one can define the objects
\[ l_{+}^I\!=\! (\rho^I \otimes id) \, r,\;\;
l_-^I \!=\! -(\rho^I \otimes id) \, {r}^{\prime},\;\;\;
r_+^{IJ} \!=\! (\rho^I\otimes\rho^J)\,r,\;\;
r_-^{IJ} \!=\!- (\rho^I\otimes\rho^J)\,{r}^{\prime},\]
where $r^{\prime}=\sigma(r)$ with $\sigma$ being the permutation in
${\cal J}\otimes{\cal J}$.

Now we can write down the classical limit of the equations (\ref{2.1}) and
(\ref{co3}):
\beq{lim}
{}[\, \lo{1}{I_0}{\pm}\,,\,\uo{2}{J}\, ] =\,\uo{2}{J} \,
r_{\pm}^{I_0 J},  \;\;\;\;
[\, \lo{1}{I_0}{\pm}\,,\,\wo{2}{J}\, ] =
-r_{\pm}^{I_0 J}\,\wo{2}{J}.
\eeq
\begin{prop}
The relations (\ref{lim}) define generating matrices corresponding to
nondeformed co- and contravariant tensor operators, that is, each row of
$U^J_0$ and each column of $W^J_0$ form tensor operators (which transform by
the representation $\rho^J$ of the algebra ${\cal J}$) in the sense of
definitions (\ref{3}) and (\ref{co1}), respectively.
\end{prop}

{\it Proof} is analogous to the one that was given for the Propositions 1
and 2. \qed

The reformulations of Propositions 3 -- 7 for the classical case are
obvious as well. It is interesting to mention that classical analogue of
the Proposition 8 is valid for the generating matrices $\wh{U}_0^J$,
$\wh{W}_0^J$ of any representation $\rho^J$ because the matrices
$r^{IJ}_{\pm}$ obviously satisfy (\ref{cr}) (with $\chi^J$ given by
(\ref{ch2}) ).

{\it Remark}. As a consequence of (\ref{lim}) we find
\[ {}[\, \lo{1}{I}{+} - \lo{1}{I}{-} \,,\,\uo{2}{J}\, ] =
\uo{2}{J}\, c^{IJ},  \;\;\;\;
{}[\, \lo{1}{I}{+} - \lo{1}{I}{-}\,,\,\wo{2}{J}\, ] =
-\,c^{IJ}\,\wo{2}{J}, \]
where $c^{IJ}$ is a representation of the "tensor
Casimir element" $C=r_+ -r_-=t_a\otimes t^a$ of the algebra $\cal J$.
\vspace*{2mm}

\red {\bf f) Tensor operators for $U_q(sl(2))$.}
\vspace*{2mm}

Let us consider the case of $U_q(sl(2))$ in detail. The generators of this
quantum algebra obey the commutation relations
\beq{7}
[\,X_+ \, ,\,X_- \,]=\frac{q^{2H} - q^{-2H}}{q-q^{-1}},
\;\;\;\; [\,H\, ,\, X_{\pm} \,]=\pm X_{\pm},
\eeq
and their Hopf algebra attributes look like
\beq{8}
\ar{cc} \de H = H \otimes 1 + 1 \otimes H, \; &
\de \,X_{\pm} = X_{\pm} \otimes q^{-H} + q^{H} \otimes X_{\pm}, \; \\
 & \\ S(H) = - H, \; & S(X_{\pm}) = - q^{\mp 1} X_{\pm}. \; \er
\eeq

The universal $R$-matrix for $U_q(sl(2))$ is given by \cite{KR}
\beq{un}
{\cal R}\,= q^{2 H\otimes H}\,\sum\limits_{n\geq 0}\,
\frac{ {\omega}^n \,q^{-\frac{n(n+1)}{2}} }{[n]!}\,
\left( q^{n\,H} X^n_{+}\right) \otimes \left( q^{-n\,H} X^n_{-}\right),
\;\;\;\; \omega=q-q^{-1}.
\eeq

The fundamental $L$-operators obtained from ({\ref{un}) by formulae
(\ref{lop}) are
\beq{2.9}
L_+ = \left(\ar{cc} q^{H} & \omega q^{-1/2} X_- \\
0 & q^{-H} \er \right),\;\;\;
L_- = \left(\ar{cc} q^{-H} & 0 \\
-\omega q^{1/2} X_+ & q^{H} \er \right).
\eeq
The same expressions (\ref{2.9}) when generators are replaced by their
matrix representations of spin $s$ (we use the standard notation for
$q$-numbers: $[x]=(q^x-q^{-x})/\omega$)
\beq{20}
\ar{c} \left(\rho^s(X_+)\right)_{mn} = \dl_{m\,n-1}\sqrt{[m][2s+1-m]},\;\;\;
\left(\rho^s(X_-)\right)_{mn} = \dl_{m\,n+1}\sqrt{[m-1][2s+2-m]}, \\    \\
\left(\rho^s(H)\right)_{mn}=\dl_{mn}(s+1-m),\;\;\;m=1,..,2s+1  ,\er
\eeq
can be regarded as $R_{\,\pm}^{\frac 12 \,s}$, i.e., the $R$-matrices
which are $(\frac{1}{2},s)$-representation of (\ref{un}). In particular,
the fundamental $R$-matrices are
\beq{16}
R_{\,+}^{\frac 12 \,\frac 12}=q^{-1/2}\left(\ar{cccc} q & 0 & 0 & 0 \\
0 & 1 & \omega & 0 \\ 0 & 0 & 1 & 0 \\ 0 & 0 & 0 & q \er \right) ,\;\;\;
R_{\,-}^{\frac 12 \,\frac 12}=q^{1/2}\left(\ar{cccc} q^{-1} & 0 & 0 & 0 \\
0 & 1 & 0 & 0 \\ 0 &-\omega & 1 & 0 \\ 0 & 0 & 0 &q^{-1} \er \right) .
\eeq

After substituting the $L$-operators (\ref{2.9}) and the corresponding
$R$-matrices $R_{\,\pm}^{\frac 12 \,s}$ into (\ref{2.1}) or (\ref{2.2})
we get the relations describing tensor operators corresponding to (half-)
integer spin $s$. More precisely, from (\ref{2.1}) or (\ref{2.2}) the
following matrix equations for (2s+1)$\times$(2s+1) matrix $U^s$ arise
\beq{22}
q^{H}\,U^s\,q^{-H}= U^s\,q^{\rho^s(H)}, \;\;\;\;X_{\pm}\,U^s\,q^{H}-
q^{\mp 1}\,q^{H}\,U^s\,X_{\pm}= U^s\,\rho^s(X_{\pm}),
\eeq
where $\rho^s(X_{\pm})$ and $\rho^s(H)$ are matrices defined in (\ref{20}).
It is easy to see that these matrix equations exactly coincide with the
definition of tensor operators for $U_q(sl(2))$
\beq{9}
\ar{c}
{} [\,H\, , \, T_m^s\, ] =  m \, T_m^s, \\  \\
(X_{\pm} \, T^s_m - q^m T^s_m \, X_{\pm} )\; q^{H} =
\sqrt{[s\pm m + 1][s \mp m]} \; T_{m\pm 1}^{\;s}, \er
\; m=-s,..,s
\eeq
which is a particular case of the definition (\ref{6b}).
Equations (\ref{22}) show that the rule (\ref{2.10}) now reads (the index $n$
labels different tensor operators)
\[  U^s_{n\,1} = T_{\;\,s}^{s\,\langle n\rangle},\;\;\;U^s_{n\,2} =
T_{\,s-1}^{s\,\langle n\rangle},\;...\;,\,
U^s_{n\,2s+1} = T_{\;-s}^{s\,\langle n\rangle} . \]

Treating in the same way the relations (\ref{co3}), we get the matrix
equations
\beq{36}
q^{H}\,W^s\,q^{-H}= q^{\rho^s(H)}\,W^s, \;\;\;\;X_{\pm}\,W^s\,q^{H}-
q^{\mp 1}\,q^{H}\,W^s\,X_{\pm}= -q^{\mp 1}\,\rho^s(X_{\pm})\,W^s,
\eeq
which are equivalent to the definition of conjugated tensor operator for
$U_q(sl(2))$
\beq{38}
[\,H\, ,\, \ov{T}_m^s\, ] = - m\,\ov{T}_m^s, \;\;\;
(X_{\pm} \, \ov{T}^s_m - q^{-m} \ov{T}^s_m \,X_{\pm} )\; q^{H} =
- \sqrt{[s\mp m+1][s \pm m]}\; \ov{T}_{m \mp 1}^s.
\eeq
Equations (\ref{36}) show that the entries of $n$-th column of the matrix
$W^s$ forms conjugated tensor operator as follows
\[ W^s_{1n} = \ov{T}_{\;\,s}^{s\,\langle n\rangle},\;\;\;W^s_{2n} =
\ov{T}_{s-1}^{s\,\langle n\rangle},\;...\;,
\, W^s_{(2s+1)n} = \ov{T}_{\;-s}^{s\,\langle n\rangle} .\]

Let us give some particular realization of the matrix $W^{\frac 12}$
constructed \footnote{This matrix differs from that obtained in \cite{BF} in
signs of the off-diagonal elements because these matrices were chosen to be
consistent with slightly different forms of $L_{\pm}$.}
in \cite{BF} (here $\ga$ is an arbitrary coefficient):
\beq{44}
W^{\frac 12} \! = \! \left(\! \ar{cc} q^{\frac 12}\,z_1^{-1} [z_1\partial_1]
\,q^{-\ga\,p} q^{\frac 12 z_2\partial_2} &
-z_2\,q^{\ga\,p} q^{-\frac 12 z_1\partial_1}\!  \\ & \\
\! q^{-\frac 12}\, z_2^{-1} [z_2\partial_2] \, q^{-\ga\,p}
q^{-\frac 12 z_1\partial_1} & z_1\,q^{\ga\,p}
q^{\frac 12 z_2\partial_2}\! \er \right)\!
\frac{1}{\sqrt{[p]}},\;\;\;p=z_1\partial_1\! +\! z_2\partial_2\! +\! 1 .
\eeq
Here the entries of $W^{\frac 12}$ are regarded as operators acting on
$D_q(z_1,z_2)$ -- the space of holomorphic functions of two complex variables
(a natural realization of the model space $\cal M$ for $U_q(sl(2))$,
see, e.g., \cite{BF} for more comments) spanned on the following orthonormal
basic vectors
\[ |j,m\rangle =\frac{z_1^{j+m}z_2^{j-m}}{\sqrt{[j+m]![j-m]!}},\;\;\;\;j=
\ar{lllll}0,&\!\!\frac{1}{2},&\!\!1,&\!\!\frac{3}{2},&\!\! ...\er \;\;\;\;
m=-j,..,j. \]
The realization of the generators (\ref{7}) on $D_q(z_1,z_2)$ is given by
\beq{45}
X_+ = z_1 \, z_2^{-1} [z_2 \partial_2 ],\;\;\;
X_- = z_2 \, z_1^{-1} [z_1\partial_1 ], \;\;\;
q^H = q^{\frac{1}{2} (z_1 \partial_1 - z_2 \partial_2)}
\eeq
and $p$ which enters (\ref{44}) is the operator of spin, i.e.,
$p\,|j,m\rangle=(2j+1)|j,m\rangle$.

One can verify that (\ref{44}) obeys the defining relations
(\ref{co3}) with $R$-matrices given by (\ref{16}) and entries of the
(\ref{2.9}) taken in the representation (\ref{45}). Therefore, the columns
of (\ref{44}) are two (nonequivalent) contravariant tensor operators.

Finally, in the limit $q\!\rar\! 1$ one gets from (\ref{44}) an example of
generating matrix $W^{\frac 12}$ for the nondeformed algebra $sl(2)$
\beq{44b}
 W^{\frac 12} \! = \! \left(\! \ar{rr} \partial_1 \, &  -z_2 \\
\! \partial_2 & z_1  \er \right)\!
\frac{1}{\sqrt{p}},\;\;\;p=z_1\partial_1\! +\! z_2\partial_2\! +\! 1 .
\eeq

\vspace*{3mm}

\red {\large\bf 3. Discussion}
\vspace*{2mm}

$\!\!\bullet $ In the present paper we gave the systematical description of
tensor operators for classical and quantum Lie algebras in the $R$-matrix
language. Of course, there exist other possible ways to formulate
the most of statements above (see e.g. \cite{BMT,Bi,MS}).

For example, alternative definitions can be given in the language developed
in \cite{MS} as follows. Let $\rho^I$ be a certain representation of
${\cal J}_q$. Define the maps
\[  \mu^I\,(\,\xi\,)=(\rho^I\otimes id)\,\de(\,\xi\,),
\;\;\;\xi\in{\cal J}_q     . \]
These objects can be regarded as matrices in auxiliary space:
$\mu^I\,(\,\xi\,)\in{\rm End}(V^I)\otimes {\cal J}_q$. Now the statement
that the set of operators $\{{\Psi}^I_j\}$ is a tensor operator looks like
\[ {\Psi}^I_m \, \xi = \mu^I_{km}\,(\,\xi\,) \, {\Psi}^I_k  .\]
The later relation is equivalent to our definition (\ref{2.1}) due to the
quasitriangularity (\ref{qt2}) of the universal $R$-matrix.
The analogue of the fusion formula (\ref{f3}) now reads
\[ {\Psi}^K_m=\left( {\rm C}[IJK]\,{\Psi}^I \otimes {\Psi}^J \right)_{m}, \]
where the Clebsch-Gordan maps were defined in (\ref{t3})-(\ref{t5}). Here
these CG-maps play the same role as the quantum projectors do in
(\ref{f1}), (\ref{f3}). \\

$\!\!\bullet$  So far we have considered, in fact, $q$-deformed objects for
generic value of the deformation parameter $q$. In the case of $q$ being a
root of unity the representation theory of quantum algebras and the theory
of $q$-tensor operators are to be modified. A possible approach (the method
of truncation) was proposed in \cite{MS}. It would be interesting to formulate
the theory of "truncated" tensor-operators in the $R$-matrix language. \\

$\!\!\bullet$ Finally, let us point out that the equations (\ref{2.1}) and
(\ref{co3}) can be complemented by some equations (either of the $R$-matrix
type \cite{FB,AF,BF} or of the "functorial" type \cite {MS}) which define
commutation relations between the entries of generating matrices. Introducing
such equations, one imposes strong additional conditions on a structure
of generating matrices.

For example, the matrix (\ref{44}) does satisfy the following $R$-matrix
equations (see \cite{BF} for details)
\beq{*}
  R_{\pm}\,\W{1}{}\,\W{2}{}\,=\,\W{2}{}\,\W{1}{}\, R_{\pm}^F ,\;\;\;\;
  R_{\pm}^F \equiv (F^{\prime})^{-1}\, R_{\pm}\, F,
\eeq
where $F$ is a certain matrix depending on the full spin $p$ (this object was
discussed in different aspects in \cite{FB,BBB}). It is worth to mention that
the fusion formulae (\ref{f3}) for the generating matrices obeying the
additional relation (\ref{*}) must involve exactly that matrix $F$ which
appears in (\ref{*}).

Let us also mention that, most probably, the additional relations of type
(\ref{*}) allow to get a full set of basic tensor operators ((\ref{44}) is
an example for $U_q(sl(2))$). More on this topic will be discussed in
\cite{BS}.

\vspace*{3mm}

\red {\large\bf Appendix.}
\vspace*{2mm}

 The following construction is due to V.Schomerus (see also \cite{MS}).
Let ${\cal P}^K\in{\cal J}_q$ be a projector on the subspace ${\cal H}_K
\subset{\cal H}$ corresponding to a certain representation $\rho^K$, i.e.,
$\rho^J({\cal P}^K)=\dl_{JK}\,\rho^K(e)$. Then the projectors introduced in
(\ref{f1})-(\ref{f21}) are given by
\[ \Po{}{IJ}{K}\,=\,(\rho^I\otimes\rho^J)\Bigl(\de({\cal P}^K)\,
{\cal P}^I\otimes{\cal P}^J \Bigr).  \]

Now using the quasitriangularity property (\ref{qt2}) and the commutativity
of ${\cal P}^K$ with all the elements of ${\cal J}_q$,
\[\ar{c} \Po{23}{IJ}{K}\,\R{13}{LJ}{+}\,\R{12}{LI}{+} =
\rho^L\otimes\rho^I\otimes\rho^J \left(e\otimes {\cal P}^I\otimes{\cal P}^J\,
(id\otimes\de)\,\Bigl( (e\otimes {\cal P}^K)\,{\cal R} \Bigr)\right) =\\  \\
=\rho^L\otimes\rho^I\otimes\rho^J \left( (id\otimes\de)\,\Bigl( {\cal R}\,
(e\otimes {\cal P}^K)\Bigr)\,e\otimes {\cal P}^I\otimes{\cal P}^J \right)=
 \R{13}{LJ}{+}\,\R{12}{LI}{+}\,\Po{23}{IJ}{K}  \er \]  \\
we obtain the desired equality (\ref{f21}).

\vspace*{3mm}

\red {\large\bf Acknowledgements.}
\vspace*{2mm}

I am grateful to V.Schomerus for constructive discussions and
carefully reading of the manuscript.
I would like to thank L.D.Faddeev, P.P.Kulish, P.Schupp and
M.A.Semenov-Tian-Shansky for useful comments.
I am thankful to Prof.~G.Mack for hospitality at II. Institute for
Theoretical Physics, University of Hamburg.

This work was supported in part by Volkswagen-Stiftung grant and by ISF grant
R2H000.


\newpage

\begin{thebibliography}{0}
\vspace*{-3mm}
\footnotesize{

\bibitem{0}{E.P.Wigner, \, {\it Group theory and its applications to the
quantum mechanics of atomic spectra}. (Academic Press, 1959); \\
E.M.Loebl (ed.), \, {\it Group theory and its applications}. (Academic
Press, 1968); \\
A.O.Barut, R.Raczka, \, {\it Theory of group representations
and applications}. (Scient. Publishers, 1977); \\
L.C.Biedenharn, J.D.Louck, \, {\it Angular momentum in quantum
physics}. Encyclopedia of mathematics and its applications.
{\bf v.8} (1981).}

\bibitem{3}{ I.N.Bernstein, I.M.Gelfand, S.I.Gelfand, \, {\it Models for
representations of Lie groups}, in Proc. of I.G.Petrovsky seminars,
{\bf 2}, 3 (1976) (in Russian);
\, Funct. analysis and its appl. {\bf 9}, 61 (1975). }

\bibitem{BMT}{ D.Buchholz, G.Mack, I.Todorov,\, In: Proc. of conf. on
quantum groups, Clausthal Zellerfeld, 1989; \\
M.Nomura, \, J. Phys. Soc. Japan {\bf 59}, 439 (1990); \\
G.Mack, V.Schomerus,\, Phys. Lett. {\bf B 267}, 207 (1991);
\, Comm. Math. Phys. {\bf 134}, 139 (1990); \\
F.Pan, \, J. Phys. {\bf A 24}, L803 (1991);\\
V.Rittenberg, M.Scheunert, \, J. Math. Phys. {\bf 33}, 436 (1992).}

\bibitem{Bi}{L.C.Biedenharn, M.Tarlini, \,Lett. Math. Phys. {\bf 20}, 271
(1990);
\\ K.Bragiel, \, Lett. Math. Phys. {\bf 21}, 181 (1991);\\
L.C.Biedenharn, M.Lohe, \, Lect. Notes Math. {\bf 1510}, 197 (1992);\\
A.U.Klimyk, \, J. Phys. {\bf A25}, 2919 (1992);\\
J.F.Cornwell, to appear in J. Math. Phys.}

\bibitem{MS}{G.Mack, V.Schomerus, \, Nucl. Phys. {\bf B 370}, 185 (1992).}

\bibitem{Dr}{V.G.Drinfeld, \, {\it Quantum groups}. In: Proc. of ICM,
Berkeley 1986 (Amer. Math. Soc. Publ., 1987); Sov. Math. Dokl. {\bf 32},
No.1 (1985). }

\bibitem{J}{M.Jimbo, \, Lett. Math. Phys. {\bf 10}, 63 (1985);
 Lett. Math. Phys. {\bf 11}, 247 (1986).}

\bibitem{FRT}{L.D.Faddeev, N.Yu.Reshetikhin, L.A.Takhtajan, \, Leningrad
Math. Journal {\bf 1}, 193 (1990).}

\bibitem{RS}{ N.Yu.Reshetikhin, M.A.Semenov-Tian-Shansky, \, Lett. Math.
Phys. {\bf 19}, 133 (1990).}

\bibitem{re}{E.K.Sklyanin,\, J.Phys. {\bf A 21}, 2375 (1988); \\
\ \ P.P.Kulish, E.K.Sklyanin,\, J.Phys. {\bf A 25}, 5963 (1992).}

\bibitem{BR}{N.Burroughs, \, Comm. Math. Phys. {\bf 127}, 109 (1990); \\
\ \ M.Rosso, \, Comm. Math. Phys. {\bf 124}, 307 (1989).}

\bibitem{KR}{A.N.Kirillov, N.Yu.Reshetikhin,\,
Adv. series in Math. Phys. {\bf 11}, 202 (World Scientific, 1990). }

\bibitem{RT}{A.N.Kirillov, N.Yu.Reshetikhin, \, Comm. Math. Phys. {\bf 134},
421 (1990); \\ \ \ N.Yu.Reshetikhin, V.G.Turaev, \, Comm. Math. Phys.
{\bf 127}, 1 (1990).}

\bibitem{KS}{P.P.Kulish, E.K.Sklyanin,\, Lect.Notes Phys. {\bf 151}, 61
(1982);\\ \ \ P.P.Kulish, N.Yu.Reshetikhin, E.K.Sklyanin,\, Lett. Math.
Phys. {\bf 5}, 393 (1981).}

\bibitem{BF}{A.G.Bytsko, L.D.Faddeev, \, {\it The $q$-analogue of model
space and CGC generating matrices}, preprint q-alg/9508022 (1995). }



\bibitem{FB}{L.D.Faddeev, \, Comm. Math. Phys. {\bf 132}, 131 (1990); \\
O.Babelon, \, Comm. Math. Phys. {\bf 139}, 619 (1991).}

\bibitem{AF}{ A.Yu. Alekseev, L.D.Faddeev, \, Comm. Math. Phys. {\bf 141},
413 (1991);\\ \ \
F.Falceto, K.Gawedzki, \, J. Geom. Phys. {\bf 11}, 251 (1993).}

\bibitem{BBB}{O.Babelon, D.Bernard, E.Billey, \,{\it A quasi-Hopf algebra
interpretation of quantum 3-$j$ and 6-$j$ symbols and difference
equations}, preprint q-alg/9511019 (1995). }

\bibitem{BS}{ A.G.Bytsko, V.Schomerus, \, to appear.}
}
\end{thebibliography}
\end{document}